\documentclass{amsart}

\usepackage{amsmath,amsfonts,amssymb,amsthm,mathrsfs}
\usepackage[
            pagebackref  = false,
            colorlinks   = true,
            pdfstartview = FitV,
            citecolor=red
           ]
           {hyperref}

\usepackage{url}
\usepackage{verbatim}
\numberwithin{equation}{section}

\newcommand\be{\begin{equation}}
\newcommand\ee{\end{equation}}

\newcommand\QQ{\bf{Q}}
\newcommand\G{\bf{G}}
\newcommand\A{{\bf{A}}}
\newcommand\RR{\bf{R}}

\newcommand\oo{\mathcal{O}}

\DeclareMathOperator{\GL}{GL}

\newtheorem{thm}{Theorem}

\theoremstyle{definition}

\newtheorem{rem}[thm]{Remark}

\title{On the balanced Vorono\"{i} formula for GL$_N$}

\author{Tian An Wong}
\address{Unversity of British Columbia, Vancouver, Canada}
\email{wongtianan@math.ubc.ca}


\begin{document}
\maketitle

\begin{abstract}
Miller and F. Zhou have proved a balanced Vorono\"i summation formula for GL$_N$ over $\QQ$, which allows one to control the dimensions of the Kloosterman sums appearing on either side of the Vorono\"i  formula. In this note, we  prove a balanced Vorono\"i formula over an arbitrary number field, starting with the Vorono\"i summation formula of A. Ichino and N. Templier over number fields, allowing one to extend recent results on spectral reciprocity laws to number fields, in special cases. 
\end{abstract}

\section{Introduction}

\subsection{}

The Vorono\"i summation formula is an equality between a weighted sum of Fourier coefficients of an automorphic form twisted by an additive character, and a dual weighted sum of Fourier coefficients of the dual form twisted by a Kloosterman sum. The Vorono\"i formula for $\GL_2$ is a basic tool in the study of automorphic forms, while more general applications have followed with the more general formulas for $\GL_N$ proved by Goldfeld-Li \cite{GL} and Miller-Schmid \cite{MS} over $\QQ$, and Ichino-Templier \cite{IT} over any number field $F$ which, importantly, removes any ramification assumptions in the previous cases.

A balanced formula on $\GL_N$ was first obtained by Zhou \cite{Zhou} under certain restrictions, then later  in general by Miller-Zhou \cite{MZ}, in which the lengths of the hyper-Kloosterman sums on either side of the formula can be chosen in a `balanced' manner. That is, the dimensions of the hyper-Kloosterman sums on either side of the summation formula can be taken to be $M'$ and $N'$ respectively, where $M'+L'+2 = N$. 

This was applied in the recent work of Blomer-Li-Miller \cite{BLM} to prove a spectral reciprocity law via a so-called `Kuznetsov-Vorono\"i-Kuznetsov triad' for a spectral sum of automorphic $L$-functions on $\GL_4\times \GL_2$ as follows: a Kuznetsov trace formula on $\GL_2$ is applied, and then the balanced Vorono\"i formula for $\GL_4$ is used on the geometric side, and the Kuznetsov formula is applied again to the dual geometric side, to give a dual spectral sum. As an application, the authors prove a non-vanishing result for automorphic $L$-functions on $\GL_4\times \GL_2$. A modified version has also been developed in Blomer-Khan \cite{Bkhan}, and is used to bound moments of twisted $L$-functions on $\GL_4$. The mechanics of the spectral reciprocity law suggest that a general formula may exist for $\GL_{2N}\times \GL_{N}$. Unfortunately, even for $N=3$ one finds that the Kuznetsov formula involves Kloosterman sums of varying lengths, which prevents a direct application of the balanced Vorono\"i formula as in the $N=2$ case.

\subsection{}

In this paper, we generalise the balanced Vorono\"i to a general number field. Besides allowing for extensions of the results on spectral reciprocity laws to number fields in special cases, another feature of our work is that rather than Kloosterman sums, more general Kloosterman integrals appear on either side of the balanced formula, which allows for the possibility of a more flexible relative trace formula, which involves Kloosterman integrals, to be used in place of the Kuznetsov trace formula. 

A second motivation for our study comes from a somewhat different source. Recent developments with regards to the conjectures of Braverman and Khazdan \cite{BK1} such as \cite{N1, N3, BNS} developing geometric methods to generalise the theory of Godement and Jacquet \cite{GJ}, which proves the functional equation of standard automorphic $L$-functions on $\GL_N$ using Poisson summation. In particular, \cite{N3} proposes a construction of the conjectural $\rho$-Fourier transform $\mathscr F^\rho$, which generalizes the Hankel transform that occurs in the Vorono\"i formula for $\GL_2$. The existence of balanced Vorono\"i formulas then suggests that the $\rho$-Poisson summation formula of the form
\[
\sum_{\gamma\in G(F)} \phi(\gamma)  = \sum_{\gamma\in G(F)}\mathscr F^\rho(\phi)(\gamma),
\]
where $\phi$ belongs to a certain $\rho$-Schwartz space $\mathscr S^\rho(G(\A_F))$, can be again `balanced' in a similar manner, and it would be interesting to explore potential applications to the analytic theory of $L$-functions. 

\subsection{Main result} Our method essentially follows that of Miller-Zhou, where instead of starting with the Vorono\"i formula of Miller-Schmid \cite{MS} over $\QQ$ we use the more general formula of Ichino-Templier \cite{IT}, and avoid the use of multiple Dirichlet series. The key observation is that the proof of the balanced Vorono\"i formula reduces to the usual Vorono\"i formula through a series of character sums, parallel to the repeated use of the crucial identity \cite[Lemma 3.2]{MZ}.  

Let $N+2=L+M$. Let $T$ be the maximal torus of diagonal matrices in GL$_N$, and $T^L, T^M$ disjoint sub-tori of dimensions $L-1$ and $M-1$ respectively, so that $T \simeq T^L\times T^M$.  The case $L=1$ then reduces to the ordinary Vorono\"i summation formula.  Then referring to Section \ref{defs} below for the definitions and notations, our balanced Vorono\"i formula is as follows.

\begin{thm}
\label{balancedvoronoi}
Let $F$ be a number field. Let $\pi=\otimes_v\pi_v$ be an irreducible cuspidal automorphic representation of $GL_N(\A_F)$, and let $S$ be the set of places of $F$ over which $\pi_v$ is ramified.  For any $\zeta\in \A^S_F$ and $\omega_S \in C_c^\infty(F_S^\times)$, we have
\begin{align}
\label{bvformula}
&\sum_{t\in T^M_\zeta/T^M_\circ} \sum_{\gamma\in F^\times} Kl_M(\gamma\zeta,t) W_\circ^{S}\Big( \begin{pmatrix}\gamma & \\ & 1_{N-1}\end{pmatrix}\Big) w_{S}(\gamma)\\
&=
\sum_{s\in T^L_\zeta/T^L_\circ} c(\zeta,s) \det(s) \sum_{\gamma\in F^\times} Kl_L(\gamma\zeta^{-1},s) \tilde{W}_\circ^{S}\Big( \begin{pmatrix}\gamma & \\ & 1_{N-1}\end{pmatrix}a(s)\Big) \tilde{w}_{S}(\gamma),\notag
\end{align}
where $c(\zeta,s)$ is defined in \eqref{const}.
\end{thm}
\noindent In principle, one should be able to apply this formula to obtain generalizations spectral reciprocity formulae, for example of \cite[Theorem 3]{BLM} in the case $N=2$ to totally real number fields, using the relevant Kuznetsov formula of Bruggeman and Miatello or relative trace formula. We will discuss this briefly at the end in Remark \ref{blm}. 

\begin{rem}
We briefly describe how the notation in \cite[Theorem 1.1]{MZ} can be compared to ours. First, for the same $N$, their parameters are chosen such that $M'+L'+2 = N$. Our choice of $M+L-2=N$ here differs from theirs due to the convention on hyper-Kloosterman sums used in \cite{IT}. Their balanced Vorono\"i formula takes the form:
\begin{align*}
&\sum_{\bf D|Q} D^{M}_1\dots D_M^1  \sum_{n=1}^\infty Kl_M(\bar a,n,c;{\bf Q}, {\bf D}) A({\bf q},{\bf D}, n)\omega\Big(\frac{nD_1^{M+1}\dots D_M^2}{q_1^L\dots q_L^1}\Big) \\ 
&= \sum_{\bf d|q}\frac{d_1^L\dots d_L^1}{c^{L+1}} \sum_{n=1}^\infty\sum_{\epsilon = \pm}  Kl_L(a,\epsilon n,c; {\bf q},{\bf d})   \tilde{A}({\bf Q},{\bf d},n)\Omega\Big(\frac{(-1)^M\epsilon n d_1^{L+1}\dots d_L^2}{c^NQ_1^M\dots Q_M^1}\Big).
\end{align*}
Letting $F=\QQ$, we specialise $\psi(x_v)$ to be  $e^{-2\pi i x_\infty}$ for $x_\infty\in \RR$, and $e^{2 \pi i x_p}$ for $x_p\in {\QQ}_p$. Our $\zeta$ corresponds to $\frac{\bar a}{c}$, and the set of places $S$ are the prime divisors of $c$. Our $\gamma\in F^\times$ correspond to the arguments of $\omega$ and $\Omega$ above.  Our $t \in T^M_\zeta/T^M_\circ$ corresponds to their sequence of positive integers $d_1,\dots,d_M$, where up to units we have $(t_2,\dots,t_{M-2},t_{M-1})$ equal to $\frac1 c(d_1d_2\dots d_{M}, \dots, d_1d_2, d_1)$, and similarly $s \in T^L_\zeta/T^L_\circ$ corresponds to $D_1,\dots, D_L$. Their hyper-Kloosterman sum $Kl_N(a,n,c;{\bf q,d})$ corresponds to $Kl_N(\gamma\zeta^{-1},t)$ as outlined in \cite[p.72]{IT}.  Finally, the Fourier coefficients $A$ correspond to $W_{\circ f}$ and $\tilde W_{\circ f}$ up to normalization as in \eqref{FC}, and our functions $\omega,\tilde\omega$ correspond to their $\omega,\Omega$ respectively, though their test function $\omega$ is compactly supported on $(0,\infty)$.
\end{rem}


\section{The Vorono\"i formula of Ichino-Templier}
\label{defs}

\subsection{}

Let $F$ be a number field, and $\A=\A_F$ the ring of adeles. Also let $F_v$ be a completion of $F$ at a prime $v$, with ring of integers $\oo_v$. Fix a non-trivial additive character $\psi=\otimes_v\psi_v$ of $F\backslash \A$. Let $\pi=\otimes_v\pi_v$ be an irreducible cuspidal automorphic representation of $GL_N(\A_F)$, $n\ge2$, and let $S$ be the set of places of $F$ over which $\pi_v$ is ramified. Let $\A^S$ be the adeles with trivial component above $S$. Define the unramified Whittaker function of $\pi^S$ to be $W_\circ^S = \prod_{v\not\in S}W_{\circ v}$, and similarly for the contragredient representation $\tilde{\pi}^S$ we write $\tilde{W}_{\circ}^S$, where 
\[
\tilde{W}_{\circ}^S(g) = W_{\circ}^S(w^tg^{-1})
\]
for all $g\in GL_N(\A^S)$, and $w$ is the long Weyl element of GL$_N$. Over $\QQ$, they are related to the Fourier coefficients $A(m_1,m_2,\dots,m_{N-1})$ of $\pi$ by the following relation:
\be
\label{FC}
\prod_{p<\infty} W_p(\Delta_m)= \prod_{i=1}^{N-1}|m_i|^{-i(n-i)/2}A(m_1,m_2,\dots,m_{N-1}),
\ee
where 
\[
\Delta_m= \text{diag}(m_1\dots m_{N-1}, m_2\dots m_{N-1},\dots, m_{N-1},1)
\]
is a diagonal matrix in $GL_N(\QQ)$.

\subsubsection{Measures}

Throughout, we make the following choices of measures. The measure $dx_v$ on the local field $F_v$ is chosen to be self-dual with respect to the fixed additive character $\psi_v$. Fix a non-zero differential form $\omega$ in Hom$_F(\wedge^\text{top}\text{Lie}(U),F)$ and also for $Y$, so that $\omega_v$ and $dx_v$ determine a measure on Lie$(U)(F_v)$, hence an invariant measure on $U(F_v)$. The product of these measures gives the Tamagawa measure.

\subsubsection{Generalised Bessel transforms}

Define for each $w_v \in C_c^\infty(F_v^\times)$ a dual function $\tilde{\omega}_v$ such that
\begin{align*}
&\int_{F^\times_v}\tilde{\omega}_v(y) \chi(y)^{-1}|y|^{s-\frac{N-1}{2}} dy \notag\\
&=  \chi(-1)^{N-1}\gamma(1-s,\pi_v\times\chi,\psi_v)\int_{F^\times_v} w_v(y)\chi(y)|y|^{1-s-\frac{N-1}{2}}dy
\end{align*}
for all Re$(s)$ large enough and any unitary character $\chi$ of $F^\times_v$. This defines $\tilde{\omega}_v$ uniquely in terms of $\pi_v,\psi_v,$ and $\omega_v$, independent of the choice of Haar measure $dy$. Note that $\tilde{\omega}_v(x)$ is smooth of rapid decay, but not necessarily compactly supported, as $|x|\to\infty$, which is important for the convergence of the dual sum.

\subsubsection{Kloosterman integrals}

Define for any $\gamma_v,\zeta_v \in F^\times_v$, the hyper-Koosterman integral,
\[
K_v(\gamma_v,\zeta_v,\tilde{W}_{\circ v}) :=
|\zeta_v|^{N-2} \int_{U_\tau^-(F_v)}\overline{\psi}_v(u_{N-2,N-1})\tilde{W}_{\circ v}(\tau u) du
\]
where
\[
\tau = \begin{pmatrix}&1& \\ 1_{N-2}&&\\ &&1 \end{pmatrix} \begin{pmatrix}1_{N-2}&&\\ &-\gamma_v\zeta_v^{-1}&\\ &&-\zeta\end{pmatrix},
\]
and set 
\[
K_R(\gamma,\zeta,\tilde{W}_{\circ R})=\prod_{v\in R}K_v(\gamma_v,\zeta_v,\tilde{W}_{\circ v})
\]
for $\gamma,\zeta \in \A^\times_R$. It relates to hyper-Kloosterman sums as follows. Let $T$ be the maximal torus of diagonal matrices in GL$_N$, then 
\[
K_v(\gamma_v,\zeta_v,\tilde{W}_{\circ v})=
|\zeta_v|^{N-2} \sum_{T(F)^+/T(\oo_v)} \tilde{W}(t)Kl_N(\gamma_v\zeta_v^{-1},t)
\]
where the sum is taken over elements $t=(t_1,\dots,t_N)$ in $T(F_v)^+/T(\oo_v)$ such that 
\[
1\leq |t_2|\leq \dots \leq |t_N|=|\zeta_v|, \text{ and }|t_1t_2\dots t_{N-1}|=|\zeta_v|.
\]
Here $Kl_N(\gamma\zeta^{-1},t)$ is the hyper-Kloosterman sum of dimension $N-1$ and can be expressed as
\[
\sum_{v_{N-1}\in t_{N-1}\oo^\times_v/\oo_v} \cdots \sum_{v_{2}\in t_2\oo^\times_v/\oo_v}
\psi(v_{N-1}+\dots + v_2)\psi((-1)^n\gamma\zeta_v^{-1}v_2^{-1}\dots v^{-1}_{N-1})
\]
by \cite[Corollary 6.7]{IT}.

\subsubsection{Vorono\"i formula}

We can now state the main result of Ichino and Templier \cite[Theorem 1]{IT}, which will be the basis for our balanced Vorono\"i formula.

\begin{thm}[Ichino-Templier]
\label{IchinoTemplier}
Let $\zeta\in \A^S_F$, and $R$ the set of places $v$ such that $|\zeta_v|>1$. Then with notation as above, we have
\begin{align*}
&\sum_{\gamma \in F^\times} \psi(\gamma\zeta)W_\circ^S\Big( \begin{pmatrix}\gamma & \\ & 1_{N-1}\end{pmatrix}\Big) w_S(\gamma)\\
&=
\sum_{\gamma \in F^\times}K_R(\gamma, \zeta,\tilde{W}_{\circ R})\tilde{W}_\circ^{R\cup S}\Big( \begin{pmatrix}\gamma & \\ & 1_{N-1}\end{pmatrix}\Big) \tilde{w}_S(\gamma),
\end{align*}
for any $\omega_S \in C_c^\infty(F_S^\times)$.
\end{thm}

\noindent From the preceding discussion, we can expand the right-hand side along the maximal torus $T$ to obtain an expression in terms of Kloosterman sums
\[
\sum_{t\in T_\zeta/T_\circ} \sum_{\gamma\in F^\times} Kl_N(\gamma\zeta^{-1},t) \tilde{W}_\circ^{S}\Big( \begin{pmatrix}\gamma & \\ & 1_{N-1}\end{pmatrix}a(t)\Big) \tilde{w}_S(\gamma),
\]
here $T_\zeta$ denotes the set of $(t_2,\dots,t_{N-1})$ in $F_R^{N-2}$ such that 
\[
1\leq |t_{2}|_v\leq \dots \leq |t_{N-1}|_v \leq |\zeta|_v
\]
for all $v\in R$. Here $T_\circ = (\oo^\times_R)^{N-2}$ and $a(t)$ is the diagonal matrix $(t_1,\dots,t_N)$ in $T(\A_R)/T(\oo_R)$ uniquely completed such that $|t_{N}|_v = |\zeta|_v$ and $|t_1\cdots t_N|_v=1$ for all $v\in R$. Taking $F=\QQ$, and $\pi$ to be unramified at every finite prime, this recovers the main result of \cite{MS} (see \cite[Theorem 2]{IT}).

\section{Proof of Theorem \ref{balancedvoronoi}}

We are now ready to prove a balanced Vorono\"i formula over an arbitrary number field, which specialises to Theorem \ref{IchinoTemplier} at $M=0$.  First, we open up the hyper-Kloosterman sum on the left-hand side of \eqref{bvformula}, and then bring in the $\gamma$ sum,
\begin{align}
\sum_{t\in T^M_\zeta/T^M_\circ} \sum_{\substack{v_{M-1}\in t_{M-1}\oo^\times_R/\oo_R\\ \dots \\ v_{2}\in t_{2}\oo^\times_R/\oo_R}}& \psi(v_{M-1}+\dots + v_2) \ \times  \notag\\
&\sum_{\gamma\in F^\times} \psi((-1)^M\gamma\zeta^{-1}v_2^{-1}\dots v^{-1}_{M-1})W_\circ^{S}\Big( \begin{pmatrix}\gamma & \\ & 1_{N-1}\end{pmatrix}\Big) w_{S}(\gamma).
\end{align}
Note that interchanging the summation is justified by the compact support of the test function $\omega_S$. Then applying the Vorono\"i summation of Theorem \ref{IchinoTemplier} to the inner sum, we obtain the dual expression
\begin{align}
\label{innervoronoi}
\sum_{t\in T^M_\zeta/T^M_\circ} \sum_{\substack{v_{M-1}\in t_{M-1}\oo^\times_R/\oo_R\\ \dots \\ v_{2}\in t_{2}\oo^\times_R/\oo_R}}& \psi(v_{M-1}+\dots + v_2) \ \times  \notag\\
&\sum_{s\in T_\zeta/T_\circ} \sum_{\gamma\in F^\times} Kl_N(\gamma\breve\zeta^{-1},s) \tilde{W}_\circ^{S}\Big( \begin{pmatrix}\gamma & \\ & 1_{N-1}\end{pmatrix}a(s)\Big) \tilde{w}_S(\gamma),
\end{align}
where we have denoted $\breve\zeta:=(-1)^M\zeta^{-1} v_2^{-1}\dots v_{M-1}^{-1}$. Recall that here $T$ is the maximal split torus in $G$, so that relabelling indices if necessary, we may decompose any $s\in T_\zeta/T_\circ$ into $s=s_1s_2$ where 
\begin{align*}
s_1&=(t_1,\dots, t_{L-1}) \in T^L_\zeta/T^L_\circ, \\
 s_2&=(t_L,\dots, t_{N-1}) \in T^M_\zeta/T^M_\circ,
\end{align*}
such that 
\begin{align*}
&1\leq |t_{1}|_v\leq \dots \leq |t_{L-1}|_v \leq |\breve\zeta|_v,\\
& 1\leq |t_{L}|_v\leq \dots \leq |t_{N-1}|_v \leq |\breve\zeta|_v
\end{align*}
for all $v\in R$. Note that $s_2$ is an $(M-2)$-tuple. 

Now on the dual side, opening up the $(N-1)$-dimensional hyper-Kloosterman sum along $t_2$, down to $(L-1)$ dimension, we have
\[
Kl_N(\gamma\breve\zeta^{-1},s_1s_2)=
\sum_{u_{N-1}\in t_{N-1}\oo^\times_R/\oo_R} \cdots \sum_{u_L\in t_L\oo^\times_R/\oo_R}
\psi(u_{N-1}+\dots + u_L)Kl_L(\gamma\zeta^{-1}_L,s_1)
\]
where 
\[
\zeta_L = (-1)^{N-L}\breve\zeta u_{L}\dots u_{N-1} = \zeta^{-1} v_2^{-1}\dots v_{M-1}^{-1}u_{L}\dots u_{N-1}.
\]
Only the innermost sum over $\gamma$ is infinite, so we may rearrange the order of summation by pairing the $v_{M-1}$ sum with the $u_L$ sum, the $v_{M-2}$ sum with the $u_{L+1}$ sum, and so on. Separating the $s_1$ and $s_2$ sums, we write \eqref{innervoronoi} as 
\begin{align}
\label{pairs}
\sum_{t\in T^M_\zeta/T^M_\circ} &\sum_{\substack{v_{M-1}\in t_{M-1}\oo^\times_R/\oo_R\\ \dots \\ v_{2}\in t_{2}\oo^\times_R/\oo_R}}
\sum_{s_2\in T^M_\zeta/T^M_\circ} \sum_{\substack{u_{N-1}\in t_{N-1}\oo^\times_R/\oo_R\\ \dots \\ u_{L}\in t_{L}\oo^\times_R/\oo_R}}
\psi(v_{M-1}+u_{L})\dots \psi(v_{2}+u_{N-1})  \notag\\
& \times\sum_{s_1\in T^L_\zeta/T^L_\circ} \sum_{\gamma\in F^\times} Kl_L(\gamma\zeta^{-1}_L,s_1) \tilde{W}_\circ^{S}\Big( \begin{pmatrix}\gamma & \\ & 1_{N-1}\end{pmatrix}a(s_1)\Big) \tilde{w}_S(\gamma),
\end{align}
where the third sum is over $s_2=(t_L,\dots,t_{N-1})$ as above. It remains then to evaluate the first line, noting that it is independent of the second line except for $\zeta_L$. To treat the first four sums, we  separate also the sum on $t$ in $T^M_\zeta/T^M_\circ$  into its $(M-2)$ components $(t_2,\dots,t_{M-1})$ such that 
$1\le |t_{1}|_v \le \dots \le |t_{M-1}|_v\le |\zeta|_v$
for all $v\in R$. We will omit the subscript on $|\cdot|_v$ when the context is clear. So the first two sums of \eqref{pairs} then reads for every fixed $s_2, u_L, u_{L+1},\dots, u_{N-1}$ as follows:
\be
\label{pair}
\sum_{\substack{ t_{M-1}\in F_R \\ |t_{M-1}|\le |\zeta_R|}} \sum_{v_{M-1}\in t_{M-1}\oo^\times_R/\oo_R} \psi(v_{M-1}+u_{L})
\dots
\sum_{\substack{ t_{2}\in F_R \\ |t_{2}|\le |t_3|}} \sum_{v_{2}\in t_{2}\oo^\times_R/\oo_R} \psi(v_{2}+u_{N-1}).
\ee
Consider then the first pair. We observe that for each fixed  $s_2\in T^M_\zeta/T^M_\circ$ and $u_{L}\in t_{L}\oo^\times_R/\oo_R$, the sum:
\[
\label{lemma}
\sum_{ |t_{M-1}|\le |\zeta_R|}
\sum_{v_{M-1}\in t_{M-1}\oo^\times_R/\oo_R}
\psi(v_{M-1}+u_{L}) \\
\]
is nonzero only if $t_{M-1} = t_L, u_L \equiv - v_{M-1}\text{ (mod } \oo_R)$. 
To see this, simply observe that
\[
\sum_{v_{M-1} \in t_{M-1} \oo_R^\times / \oo_R } \psi(v_{M-1}+u_L) = |t_{M-1}|_R
\]
if $t_{M-1} = t_L$ and $u_L \equiv - v_{M-1}$ modulo $\oo_R$, and is zero otherwise. We note that this is the analogue of Lemma 3.2 of \cite{MZ}. This implies that $ v_{M-1}^{-1}u_{L}\equiv-1$ mod $\oo_R$ in $\zeta_L$. Moving on to the second pair, for fixed $t_{M-1}$ and $u_{L+1}$, 
\[
\sum_{ |t_{M-2}|\le |t_{M-1}|} \sum_{v_{M-2}\in t_{M-2}\oo^\times_R/\oo_R} \psi(v_{M-2}+u_{L+1}) 
\]
we see that the sum is again nonzero only if $t_{M-2} = t_{L+1}$ and $u_{L+1}\equiv - v_{M-2}$ mod $\oo_R$, and zero otherwise. Applying this $M-2$ times, we collect the evaluated sums \eqref{pairs} into a constant equal to the product of 
\[
|t_{M-1} \dots t_2|_R = \det(s_1)
\]
and the cardinality 
\be
\label{const}
c(\zeta,s_1)= \#\{t\in T^M_\zeta/T^M_\circ: 1\le |t_{1}|_v \le \dots \le |t_{M-1}|_v\le |\zeta|_v, v\in R\}.
\ee
 And finally the sum reduces to 
\[
\sum_{s_1\in T^L_\zeta/T^L_\circ} c(\zeta,s_1) \det(s_1)\sum_{\gamma\in F^\times} Kl_L(\gamma\zeta_L^{-1},s_1) \tilde{W}_\circ^{S}\Big( \begin{pmatrix}\gamma & \\ & 1_{N-1}\end{pmatrix}a(s_1)\Big) \tilde{w}_{S}(\gamma)
\]
as desired. 

\begin{rem}
\label{blm}
The spectral reciprocity formula obtained in \cite{BLM} corresponds to the case $N=2$ and $L=M=3$. Our balanced formula in Theorem \ref{bvformula} then consists of Kloosterman sums on either side, attached to tori $T^L$ and $T^M$ in GL$_3$ over $F$. Fixing a cuspidal automorphic representation $\Pi$ on GL$_4(\A)$ and a suitable test function $h$, we consider the spectral mean value roughly of the form
\[
\sum_{\pi}\frac{L(\frac12,\Pi\times \pi)}{L(1,\pi, \text{Sym}^2)}h(t_\pi) +\frac{1}{2\pi}\int_{-\infty}^\infty\frac{L(\frac12+it,\Pi)L(\frac12-it,\Pi)}{|\zeta(1+2it)|^2}h(t) dt
\]
where $\pi$ runs over cuspidal automorphic representations of GL$_2(\A)$. (In \cite{BLM} the sum is twisted by the parity of the eigenvalue of $\pi_\infty$ under the reflection $z\mapsto -\bar{z}$.) We expand the $L$-functions as usual to obtain 
\[
L(s,\Pi\times \pi) = \sum_{n,m>1}\frac{a_\Pi(n,m)a_\pi(n)}{n^sm^{2s}}
\]
and similarly for $L(s+it,\Pi)L(s-it,\Pi)$, and applying known convexity the original expression is seen to be absolutely convergent. 

To apply the so-called Kuznetsov-Voronoi-Kuznetsov triad over $F$, one may use either the relative or Kuznetsov trace formula. In the former case, one has Kloosterman integrals as in \cite[\S3.3]{KL} roughly of the form
\[
I_\delta(h) = \int_{H_\delta (F)\backslash H(\A)} h(n_1^{-1}\delta n_2)\overline{\theta_{m_1}(n_1)}\theta_{m_2}(n_2) dn_1dn_2
\]
where $dn_1$ and $dn_2$  are quotient measures obtained from $H(\A) = N(\A) \times N(\A)$, $H_\delta$ is the stabilizer of $\delta\in N(F) \backslash \mathrm{GL}_2(F)/N(F)Z(F)$, and $\theta_{m_1}$ and $\theta_{m_1}$ are additive characters on $N(F)\backslash N(\A)$. (See also \cite[\S7]{KL2}.)Up to a constant, this is equal to
\[
\sum_{c\in N{\bf Z}^+}\frac{S(m_1,m_2;\mathfrak n,c)}{c}J_{k-1}\left(\frac{4\pi \sqrt{\mathfrak n m_1m_2}}{c}\right).
\]
If we decompose the adelic integral into
\[
I_\delta(h) = I_{\delta_S}(h) I_{\delta^S}(h),
\]
then we see that $I_{\delta_S}(h)$, which specializes to $J_{k-1}$, corresponds to our generalized Bessel transform $\tilde{\omega}_S(\gamma)$, whereas $I_{\delta_S}(h)$ which specializes to the Kloosterman sum corresponds to the Kloosterman integral $Kl_3(\gamma\zeta,t)$. To execute the spectral reciprocity will require more careful estimates as in \cite{BLM}, and a proper generalization of the parity condition which allows the Kuznetsov formula to be inverted as Motohashi's work is over ${\bf Q}$ \cite[Theorems 2.3 and 2.5]{Mo}. 
\end{rem}

\noindent\emph{Acknowledgments.} The author thanks Giacomo Cherubini for helpful discussions and comments on a preliminary version of this paper, and the referee for comments improving the content and exposition of the paper.

\bibliographystyle{alpha}
\bibliography{specrec_90_bib}

\end{document}